\pdfoutput=1
\documentclass[11pt,a4paper]{amsart}%
\usepackage{hyperref}
\usepackage{amsfonts}
\usepackage{amsmath}
\usepackage{amssymb}
\usepackage{amsthm}
\usepackage{graphicx}
\usepackage{algorithm}
\usepackage[noend]{algorithmic}
\usepackage{float}
\usepackage{color}%
\theoremstyle{plain}

\theoremstyle{definition}

\theoremstyle{remark}

\numberwithin{equation}{section}
\theoremstyle{plain}

\DeclareMathOperator{\reg}{reg}
\DeclareMathOperator{\codim}{codim}
\DeclareMathOperator{\hilb}{Hilb}
\DeclareMathOperator{\depth}{depth}

\begin{document}
\title[Decomposition of Monomial Algebras: Applications and Algorithms]{Decomposition of Monomial Algebras:\linebreak Applications and Algorithms}
\author{Janko B\"{o}hm}
\address{Fachbereich Mathematik, TU Kaiserslautern,
67663 Kaiserslautern, Germany}
\email{boehm@mathematik.uni-kl.de}
\author{David Eisenbud}
\address{Department of Mathematics, University of California, Berkeley, CA 94720, USA}
\email{de@msri.org}
\author{Max J. Nitsche}
\address{Max\mbox{\;}Planck\mbox{\;}Institute\mbox{\;}for\mbox{\;}Mathematics\mbox{\;}in\mbox{\;}the\mbox{\;}Sciences,\mbox{\;}Inselstrasse\mbox{\;}22,\mbox{\;}04103\mbox{\;}Leipzig,\mbox{\;}Germany}
\email{nitsche@mis.mpg.de}
\thanks{}
\date{\today}

\begin{abstract}
Considering finite extensions $K[A]\subseteq K[B]$ of positive affine
semigroup rings over a field $K$ we have developed in \cite{BENAX} an
algorithm to decompose $K[B]$ as a direct sum of monomial ideals in $K[A]$. By
computing the regularity of homogeneous semigroup rings from the decomposition
we have confirmed the Eisenbud-Goto conjecture in a range of new cases not
tractable by standard methods. Here we first illustrate this technique and its
implementation in our \textsc{Macaulay2} package \textsc{MonomialAlgebras} by
computing the decomposition and the regularity step by step for an explicit
example. We then focus on ring-theoretic properties of simplicial semigroup
rings. From the characterizations given in \cite{BENAX} we develop and prove
explicit algorithms testing properties like Buchsbaum, Cohen-Macaulay,
Gorenstein, normal, and seminormal, all of which imply the Eisenbud-Goto
conjecture. All algorithms are implemented in our \textsc{Macaulay2} package.

\end{abstract}
\maketitle

\section{Introduction}

Let $B$ be a positive affine semigroup, that is, $B$ is a finitely generated
subsemigroup of $\mathbb{N}^{m}$ for some $m$. Let $K$ be a field and $K[B]$
the affine semigroup ring associated to $B$, which can be identified with the
subring of $K[t_{1},\ldots,t_{m}]$ generated by monomials $t^{u}:=t_{1}%
^{u_{1}}\cdot\ldots\cdot t_{m}^{u_{m}}$, where $u=(u_{1},\ldots,u_{m})\in B$.
Denote by $C(B)$ and by $G(B)$ the cone and the group generated by $B$. From
now on let $A\subseteq B$ be positive affine semigroups with $C(A)=C(B)$. We
will now discuss the decomposition of $K[B]$ into a direct sum of monomial
ideals in $K[A]$. Observe that
\[
K[B]=\bigoplus\nolimits_{g\in G}K\cdot\left\{  t^{b}\mid b\in B\cap g\right\}
,
\]
where $G:=G(B)/G(A)$. Note that $C(A)=C(B)$ if and only if $K[B]$ is a
finitely generated $K[A]$-module. From this it follows that $G$ is finite, and we can compute
the above decomposition since all summands are finitely generated. Moreover,
there are shifts $h_{g}\in G(B)$ such that
\[
I_{g}:=K\cdot\left\{  t^{b-h_{g}}\mid b\in B\cap g\right\}
\]
is a monomial ideal in $K[A]$. Thus,
\[
K[B]\cong\bigoplus\nolimits_{g\in G}I_{g}(-h_{g})
\]
as $\mathbb{Z}^{m}$-graded $K[A]$-modules (with $\deg t^{b}=b$). A detailed
formulation of the algorithm computing the ideals $I_{g}$ and shifts $h_{g}$
and a more general version of the decomposition in the setup of cancellative
abelian semigroup rings over an integral domain can be found in
\cite[Algorithm 1, Theorem 2.1]{BENAX}.

Our original motivation for developing this decomposition was to provide a
fast algorithm to compute the Castelnuovo-Mumford regularity $\reg K[B]$ of a
homogeneous semigroup ring in order to test the Eisenbud-Goto conjecture
\cite{EG}. Recall that the \emph{Castelnuovo-Mumford regularity} $\reg M$ of a
finitely generated graded module $M$ over a standard graded polynomial ring
$R=K[x_{1},\ldots,x_{n}]$ is defined as the smallest integer $m$ such that
every $j$-th syzygy module of $M$ is generated by elements of degree $\leq
m+j$. Moreover, $B$ is called a \emph{homogeneous semigroup} if there exists a
group homomorphism $\deg:G(B)\rightarrow\mathbb{Z}$ with $\deg b_{i}=1$ for
$i=1,\ldots,n$, where $\hilb(B)=\{b_{1},\ldots,b_{n}\}$ is the minimal
generating set of $B$; by $\reg K[B]$ we mean its regularity with respect to
the $R$-module structure which is given by the $K$-algebra homomorphism
$R\twoheadrightarrow K[B]$, $x_{i}\mapsto t^{b_{i}}$.

The \emph{toric Eisenbud-Goto conjecture} can be formulated as follows: let
$K$ be a field and $B$ a homogeneous semigroup, then $\reg K[B]\leq\deg
K[B]-\codim K[B]$, where $\deg K[B]$ denotes the degree and $\codim
K[B]:=\dim_{K}K[B]_{1}-\dim K[B]$ the codimension. Even this special case of
the Eisenbud-Goto conjecture is largely open, for references on known results
see \cite[Section 4]{BENAX}. The regularity of $K[B]$ is usually computed from
a minimal graded free resolution. If $n$ is large this computation is very
expensive, and hence it is impossible to test the conjecture systematically in
high codimension using this method. However, choosing $A$ to be generated by
minimal generators $e_{1},\ldots,e_{d}$ of $C(B)$ of degree $1$ the regularity
can be computed as
\[
\reg K[B]=\max\{\reg I_{g}+\deg h_{g}\mid g\in G\},
\]
where $\reg I_{g}$ denotes the regularity of $I_{g}$ with respect to the
canonical $T=K[x_{1},\ldots,x_{d}]$-module structure given by
$T\twoheadrightarrow K[A]$, $x_{i}\mapsto t^{e_{i}}$. Since the free
resolution of every ideal $I_{g}$ appearing has length at most $d-1$, this
computation is typically much faster than the traditional approaches. This
enabled us to test the conjecture for a large class of homogeneous semigroup
rings by using our regularity algorithm. See \cite[Section 4]{BENAX} for details.

In Section \ref{secdec}, we illustrate, step by step, decomposition and
regularity computation for an explicit example using our \textsc{Macaulay2}
\cite{M2} package \textsc{MonomialAlgebras} \cite{BEN}. We say that 
$K[B]$ is a \emph{simplicial semigroup ring} if the cone
$C(B)$ is simplicial. In Section
\ref{simplicial}, we focus on simplicial semigroup rings $K[B]$.
  Based on the characterizations of
ring-theoretic properties given in \cite[Proposition~3.1]{BENAX} we develop
explicit algorithms for testing whether $K[B]$ is Buchsbaum, Cohen-Macaulay,
Gorenstein, seminormal, or normal. We also discuss that, by known results, all
these ring-theoretic properties imply the Eisenbud-Goto conjecture. The
algorithms mentioned are implemented in our \textsc{Macaulay2} package.

\section{Decomposition and regularity\label{secdec}}

Our \textsc{Macaulay2} package can be loaded by

\texttt{Macaulay2, version 1.4}

\texttt{with packages: ConwayPolynomials, Elimination, IntegralClosure,
LLLBases,}

\texttt{PrimaryDecomposition, ReesAlgebra, TangentCone}

\texttt{i1 : needsPackage "MonomialAlgebras"};

We discuss the decomposition at the example of the homogeneous
semigroup $B\subset\mathbb{N}^{3}$ specified by a list of generators

\texttt{i2 : B =
\{\{4,0,0\},\{2,2,0\},\{2,0,2\},\{0,2,2\},\{0,3,1\},\{3,1,0\},\{1,1,2\}\};}

\noindent As an input for our algorithm we encode this data in a multigraded
polynomial ring

\texttt{i3 : K = ZZ/101;}

\texttt{i4 : S = K[x\_1 .. x\_7, Degrees=%
$>$%
B];}

\noindent The command

\texttt{i5 : dc = decomposeMonomialAlgebra S}

\texttt{
\begin{tabular}
[c]{lllll}%
\hspace{-0.15in}\texttt{o5 = \hspace{0cm} HashTable\{} & $\mathtt{(-1,1,0)}$ &
\texttt{=$>$\{ }ideal (x$_{1},\mathtt{x}_{3}$), & $(-1,1,0)$ & \}\\
& $\mathtt{0}$ & \texttt{=$>$\{ }ideal 1, & $0$ & \}~\}
\end{tabular}
}

\noindent decomposes $K[B]$ over $K[A]$ where $A\subseteq B$ is generated by
minimal generators of $C(B)$ with minimal coordinate sum; so in the example
$A=\langle(4,0,0),(2,2,0),(2,0,2),(0,2,2),(0,3,1)\rangle$. The keys of the
hash table represent the elements of $G$ and the values are the tuples
$(I_{g},h_{g})$, hence%
\begin{equation}
K[B]\cong\left\langle \overline{x_{1}},\overline{x_{3}}\right\rangle
(-(-1,1,0))\oplus K[A] \label{iso11}%
\end{equation}
as $\mathbb{Z}^{3}$-graded $K[A]$-modules; here we write $K[A]\cong T/J$ with
$T=K[x_{1},x_{2},x_{3},x_{4},x_{5}]$ and $\overline{x_{i}}$ for the class of
$x_{i}$. Note that the on-screen output of \textsc{Macaulay2} does not
distinguish between the class and the representative. To compute $\reg I_{g}$
we will consider the standard grading on $T$:

\texttt{i6 : KA = ring first first values dc;}

\texttt{i7 : T = newRing(ring ideal KA,Degrees=%
$>$%
\{5:1\});}

\texttt{i8 : J = sub(ideal KA,T);}

\texttt{o8 : ideal of T}

\texttt{i9 : betti res J}

\texttt{o9 =
\begin{tabular}
[c]{llll}
& 0 & 1 & 2\\
total: & 1 & 3 & 2\\
0 & 1 & . & .\\
1 & . & 1 & .\\
2 & . & 2 & 2
\end{tabular}
}

\texttt{o9 : BettiTally}

\noindent The usual approach would be to obtain $\reg K[B]$ from a minimal
graded free resolution of the toric ideal $I_{B}$ with respect to the standard grading.

\texttt{i10 : IB = monomialAlgebraIdeal S;}

\texttt{o10 : ideal of S}

\texttt{i11 : R = newRing(ring IB,Degrees=%
$>$%
\{7:1\});}

\texttt{i12 : betti res sub(IB,R)}

\texttt{o12 =
\begin{tabular}
[c]{lllllll}
& 0 & 1 & 2 & 3 & 4 & 5\\
total: & 1 & 8 & 15 & 13 & 6 & 1\\
0 & 1 & . & . & . & . & .\\
1 & . & 6 & 8 & 3 & . & .\\
2 & . & 2 & 3 & . & . & .\\
3 & . & . & 4 & 10 & 6 & 1
\end{tabular}
}

\texttt{o12 : BettiTally}

\noindent We observe that $\reg K[B]=3$. With $\deg u =(u_{1}+u_{2}+u_{3})/4$
by Equation~(\ref{iso11}) it holds%
\[
\reg K[B]=\max\left\{  \reg\left\langle \overline{x_{1}},\overline{x_{3}%
}\right\rangle + 0 ,\text{ }\reg K[A] + 0\right\}  \text{.}%
\]

\noindent By \texttt{o9} we have $\reg K[A]=2$. To compute $\reg\left\langle
\overline{x_{1}},\overline{x_{3}}\right\rangle $ do:

\texttt{i13 : I1 = first (values dc)\#0}

\texttt{o13 = \hspace{0cm} ideal (x$_{1}$,x$_{3}$)}

\texttt{o13 : ideal of KA}

\texttt{i14 : g = matrix entries sub(gens I1, T);}

\texttt{\hspace{3cm}1}\texttt{\hspace{1.25cm}2}

\texttt{o15 : Matrix T <--- T}

\texttt{i15 : betti res image map(coker gens J, source g, g)}

\texttt{o15 =
\begin{tabular}
[c]{lllll}
& 0 & 1 & 2 & 3\\
total: & 2 & 5 & 4 & 1\\
1 & 2 & 2 & . & .\\
2 & . & . & . & .\\
3 & . & 3 & 4 & 1
\end{tabular}
}

\noindent Hence $\reg\left\langle \overline{x_{1}},\overline{x_{3}%
}\right\rangle =3$ and therefore $\reg K[B]=3$. Observe, that the resolution
of $K[B]$ has length $5$, whereas the ideals $I_{g}$ have resolutions of
length at most $3$. The command

\texttt{i16 : regularityMA S}

\texttt{o16 : \{3,  \{\{ideal (x$_{1},$x$_{3}$), $(-1,1,0)$\}\}\}}

\noindent provides an implementation of this approach, also returning the
tuples $(I_{g},h_{g})$ where the maximum is achieved. By
\cite[Proposition~4.1]{BENAX} we have $\deg K[B]=\#G\cdot\deg K[A]=10$ since

\texttt{i17 : degree J}

\texttt{o17 = \hspace{0cm} 5}

\noindent Moreover, $\codim K[B]=4$ since $\dim K[B]=\dim C(B)=3$. Hence the
ring $K[B]$ satisfies the Eisenbud-Goto bound.

\section{Algorithms for ring theoretic properties\label{simplicial}}

In this section we focus on simplicial semigroup rings $K[B]$. Based on the
characterizations given in \cite[Proposition~3.1]{BENAX} we develop and prove
explicit algorithms for testing whether $K[B]$ is Buchsbaum, Cohen-Macaulay,
Gorenstein, seminormal, or normal. Note that, in the simplicial case, all
these properties are independent of $K$, and they imply the Eisenbud-Goto
conjecture by results of \cite{SVBBEG2, EGCM, MNSNAX}. As an example, consider
the following homogeneous simplicial semigroup $B\subset\mathbb{N}^{3}$
specified by the generators

\texttt{i18 : B =
\{\{4,0,0\},\{0,4,0\},\{0,0,4\},\{1,0,3\},\{0,2,2\},\{3,0,1\},\{1,2,1\}\};}

\noindent We compute the decomposition of $K[B]$ over $K[A]$, where
$A=\langle(4,0,0),(0,4,0),(0,0,4)\rangle\subset B$ is generated again by
minimal generators of $C(B)$ with minimal coordinate sum.

\texttt{i19 : S = K[x\_1 .. x\_7, Degrees=%
$>$%
B]};

\texttt{i20 : decomposeMonomialAlgebra S}

\texttt{
\begin{tabular}
[c]{lllll}%
\hspace{-0.15in}\texttt{o20 : HashTable\{} & $(-1,0,1)$ & \texttt{=$>$\{
}ideal 1, & $(3,0,1)$ & \}\\
& $(-1,2,-1)$ & \texttt{=$>$\{ }ideal 1, & $(3,2,3)$ & \}\\
& $(0,2,2)$ & \texttt{=$>$\{ }ideal 1, & $(0,2,2)$ & \}\\
& $(1,0,-1)$ & \texttt{=$>$\{ }ideal 1, & $(1,0,3)$ & \}\\
& $(1,2,1)$ & \texttt{=$>$\{ }ideal 1, & $(1,2,1)$ & \}\\
& $(2,0,2)$ & \texttt{=$>$\{ }ideal (x$_{3},$x$_{1},$x$_{2}$), & $(2,0,2)$ &
\}\\
& $(2,2,0)$ & \texttt{=$>$\{ }ideal 1, & $(2,2,4)$ & \}\\
& $0$ & \texttt{=$>$\{ }ideal 1, & $0$ & \}~\}
\end{tabular}
}

\noindent Hence
\[
K[B]\cong K[A]\oplus K[A](-1)^{4}\oplus K[A](-2)^{2}\oplus\left\langle
x_{1},x_{2},x_{3}\right\rangle (-1)
\]
with respect to the standard grading induced by $\deg u =(u_{1}+u_{2}%
+u_{3})/4$. It follows that $\depth K[B]=1$, thus, $K[B]$ is not
Cohen-Macaulay. Hence $K[B]$ is also not normal by \cite{MHCM}. We can test
seminormality via the following algorithm:

\begin{algorithm}[H]
\caption{Seminormality test}
\label{alg-seminormal}
\begin{algorithmic}[1]
\REQUIRE A simplicial semigroup $B\subseteq\mathbb{N}^{m}$.
\ENSURE \texttt{true} if $K[B]$ is seminormal, \texttt{false} otherwise.
\STATE Let $e_{1},\ldots,e_{d}\in B$ be minimal generators of $C(B)$ with minimal
coordinate sum, and set $A:=\left\langle e_{1},\ldots,e_{d}\right\rangle$.
\STATE Compute $B_{A}:=\{x\in B\mid x\notin B+(A\setminus\{0\})\}$ as described
in \cite[Algorithm ~1, Step~1]{BENAX}.
\FORALL {$x\in B_{A}$}
\STATE Solve the linear system of equations $\sum_{i=1}^{d}\lambda_{i}e_{i}=x$ for $\lambda=(\lambda_{1},\ldots,\lambda_{d})\in\mathbb{Q}^{d}$.
\STATE \textbf{if} {$\left\Vert \lambda\right\Vert _{\infty}>1$}  \textbf{then return} \texttt{false}
\ENDFOR
\RETURN \texttt{true}
\end{algorithmic}
\end{algorithm}

Here, by $\left\Vert -\right\Vert _{\infty}$ we denote the maximum norm. Note
that all $\lambda_{i}$ are non-negative since $C(B)$ is a simplicial cone.
Verifying in Step $5$ the condition $\left\Vert \lambda\right\Vert _{\infty
}\geq1$ instead results in an algorithm which tests normality. Using our
package we observe that $K[B]$ is not seminormal:

\texttt{i21 : isSeminormalMA B}

\texttt{o21 : false}

\noindent The Buchsbaum property can be tested by the following algorithm. We
denote by $K[A]_{+}$ the homogeneous maximal ideal of $K[A]$.

\begin{algorithm}[H]
\begin{algorithmic}[1]
\caption{Buchsbaum test}
\REQUIRE A simplicial semigroup $B=\left\langle b_{1},\ldots,b_{n}%
\right\rangle \subseteq\mathbb{N}^{m}$.
\ENSURE \texttt{true} if $K[B]$ is Buchsbaum, \texttt{false} otherwise.
\STATE Let $e_{1},\ldots,e_{d}\in B$ be minimal generators of $C(B)$ with minimal
coordinate sum, and set $A:=\left\langle e_{1},\ldots,e_{d}\right\rangle $.
\STATE Using the (minimal) generators $e_{1},\ldots,e_{d}$ of $A$ decompose%
\[
K[B]\cong\bigoplus\nolimits_{g\in G}I_{g}(-h_{g}),
\]
where $I_{g}\subseteq K[A]$, $h_{g}\in G(B)$ and $G=G(B)/G(A)$ by \cite[Algorithm~1]{BENAX}.
\STATE \textbf{if} {$\exists g\in G$ with $I_{g}\neq K[A]$ and $I_{g}\neq K[A]_{+}$}  \textbf{then return} \texttt{false}
\STATE $H:=\left\{  h_{g}\mid g\in G\text{ with }I_{g}= K[A]_{+}\right\}  $
\STATE $C:=\left\{  b_{1},\ldots,b_{n}\right\}  \backslash\left\{0,e_{1},\ldots,e_{d}\right\}  $
\STATE $H+C:=\left\{  h_{g}+b_{i}\mid h_{g}\in H\text{, }b_{i}\in C\right\}  $
\RETURN \texttt{true} if $\left(  H+C\right)  \cap H=\emptyset$ and
\texttt{false} otherwise.
\end{algorithmic}
\end{algorithm}

\begin{proof}
By \cite[Proposition~3.1]{BENAX} the ring $K[B]$ is Buchsbaum iff each ideal
$I_{g}$ is either equal to $K[A]$, or to $K[A]_{+}$ and $h_{g}+b\in B$ for all
$b\in\hilb(B)$. So, Step~3 is correct and we may now assume that $I_{g}=K[A]$
or $I_{g}=K[A]_{+}$ for all $g\in G$. Recall that $I_{g}=\{t^{v-h_{g}}\mid
v\in\Gamma_{g}\}K[A]$, where $\Gamma_{g}=\{x\in B_{A}\mid x\in g\}$. Moreover,
note that $\{t^{v-h_{g}}\mid v\in\Gamma_{g}\}$ is always a minimal generating
set of $I_{g}$ and $h_{g}=\sum\nolimits_{k=1}^{d}\min\left\{  \lambda_{k}%
^{v}\mid v\in\Gamma_{g}\right\}  e_{k}$, where $v=\sum_{k=1}^{d}\lambda
_{k}^{v}e_{k}$ with $\lambda_{k}^{v}\in\mathbb{Q}$.

Since $h+e_{k}\in B_{A}$ for all $h\in H$ and all $k=1,\ldots,d$, we have
$H\cap B=\emptyset$. In case that $(H+C)\cap H\not =\emptyset$ we obtain
$h+b\notin B$ for some $h\in H$ and some $b\in B\setminus\{0\}$, that is,
$h+\hilb(B) \not \subseteq B$. Hence $K[B]$ is not Buchsbaum.

In case that $K[B]$ is not Buchsbaum, there is an $h\in H$ and some
$b\in\hilb(B)$ such that $h+b\notin B$. It is now sufficient to show that
$b\in C$ and $h+b\in H$. By the above argument, $b\in C$. Let $m_{k}%
=h+b+e_{k}$ for $k=1,\ldots,d$. Suppose that $m_{i}\notin B_{A}$ for some
$i\in\{1,\ldots,d\}$. Since $m_{k}-e_{k}\notin B$ for all $k=1,\ldots,d$,
necessarily $m_{i}-e_{j}\in B$ for some $j\not =i$. Consider $y=m_{j}%
-\sum_{k=1}^{d}n_{k}e_{k}\in B$ with $n_{k}\in\mathbb{N}$ such that
$\sum_{k=1}^{d}n_{k}$ is maximal. By construction $y\in B_{A}$, moreover,
$n_{j}=0$ since $m_{j}-e_{j}\notin B$. In the same way if $x=m_{i}-e_{j}%
-\sum_{k=1}^{d}n_{k}e_{k}\in B$ with $\sum_{k=1}^{d}n_{k}$ maximal, then $x\in
B_{A}$. Since $m_{i},m_{j}\in g$ for some $g\in G$, we also have $x,y\in g$.
Since $e_{1},\ldots,e_{d}$ are linearly independent, we have $\lambda_{j}%
^{y}-\lambda_{j}^{x}\geq2$. Moreover, since $t^{y-h_{g}},t^{x-h_{g}}\in K[A]$
we get that $t^{y-h_{g}}$ is not a linear form. Hence $I_{g}\neq K[A]$ and
$I_{g}\neq K[A]_{+}$, thus, $m_{k}\in B_{A}$ for all $k=1,\ldots,d$. We have
$\#\Gamma_{g}\in\{1,d\}$ by minimality, hence $\Gamma_{g}=\{m_{1},\ldots
,m_{d}\}$. By construction, $h_{g}=h+b$ and $I_{g}=K[A]_{+}$, therefore
$h+b\in H$.
\end{proof}

Note that in Step $2$ the shifts $h_{g}$ and hence the ideals $I_{g}$ are
uniquely determined since $e_{1},\ldots,e_{d}$ are linearly independent. This
is not true for arbitrary generating sets. By

\texttt{i22 : isBuchsbaumMA B}

\texttt{o22 : true}

\noindent we conclude that $K[B]$ satisfies the Eisenbud-Goto conjecture by
\cite{SVBBEG2}. Note that we can read off from the decomposition the
regularity and the Eisenbud-Goto bound: we have $\reg K[A]=0$ and
$\reg\left\langle x_{1},x_{2},x_{3}\right\rangle =1$, therefore $\reg
K[B]=\max\{0,1,2,1+1\}=2$. Moreover, $\deg K[B]$ is the number of ideals which
occur in the decomposition, hence $\deg K[B]-\codim K[B]=8-4=4$.

Note that, in case $B$ is Buchsbaum, the regularity of $K[B]$ is independent
of the field $K$ since all ideals in the decomposition are equal to the
homogeneous maximal ideal or to $K[A]$.

We finish this section by providing an algorithm for testing the Gorenstein property.

\begin{algorithm}[H]
\caption{Gorenstein test}
\label{alg-Gorenstein}
\begin{algorithmic}[1]
\REQUIRE A simplicial semigroup $B\subseteq\mathbb{N}^{m}$.
\ENSURE \texttt{true} if $K[B]$ is Gorenstein, \texttt{false} otherwise.
\STATE Let $e_{1},\ldots,e_{d}\in B$ be minimal generators of $C(B)$ with minimal
coordinate sum, and set $A:=\left\langle e_{1},\ldots,e_{d}\right\rangle $.
\STATE Using the (minimal) generators $e_{1},\ldots,e_{d}$ of $A$ decompose%
\[
K[B]\cong\bigoplus\nolimits_{g\in G}I_{g}(-h_{g}),
\]
where $I_{g}\subseteq K[A]$, $h_{g}\in G(B)$ and $G=G(B)/G(A)$ by
\cite[Algorithm $1$]{BENAX}.
\STATE \textbf{if} {$\exists g\in G$ with $I_{g}\neq K[A]$} \textbf{then return} \texttt{false}
\STATE $H:=\left\{  h_{g}\mid g\in G\right\}  $
\STATE \textbf{if} {$h\in H$ with maximal coordinate sum is not unique}  \textbf{then return} \texttt{false}
\STATE Let $h\in H$ with maximal coordinate sum.
\WHILE {$H\neq\emptyset$}
\STATE Let $h_{g}\in H$
\STATE \textbf{if} {$h-h_{g}\notin H$}  \textbf{then return} \texttt{false}
\STATE $H:=H\backslash\left\{  h_{g},\text{ }h-h_{g}\right\}  $
\ENDWHILE
\RETURN \texttt{true}
\end{algorithmic}
\end{algorithm}

\begin{proof}
By \cite[Proposition~3.1]{BENAX} the ring $K[B]$ is Gorenstein iff
$I_{g}=K[A]$ for all $g\in G$ and $H$ has a unique maximal element with
respect to $\leq$ given by $x\leq y$ if there is a $z\in B$ such that $x+z=y$.
Note that $H=B_{A}$ since $I_{g}=K[A]$ for all $g\in G$. If there is a maximal
element $h\in H$, then this element has maximal coordinate sum. If $H$ has
more than one element with maximal coordinate sum, then $H$ does not have a
unique maximal element. To complete the proof we need to show that an element
$h_{g}\in H$ satisfies $h_{g}\leq h$ iff $h-h_{g}\in H$. But this follows from
the fact that if $x\notin B_{A}$ then $x+y\notin B_{A}$ for all $x,y\in B$.
\end{proof}

Note that performing Steps $1$--$3$ of Algorithm \ref{alg-Gorenstein} (and
returning \texttt{true} afterwards) gives a test for the Cohen-Macaulay property.

\enlargethispage{6mm}
\end{document}